\documentclass[11pt]{amsart}

\usepackage[colorlinks=true, pdfstartview=FitV, linkcolor=blue,
citecolor=blue]{hyperref}

\usepackage{amssymb,amsmath,amscd}
\usepackage{graphicx}
\usepackage{a4wide}
\usepackage{tikz}
\usetikzlibrary{graphs,arrows,shapes}
\tikzset{every picture/.style={line width=0.75pt}}
\newtheorem{theorem}{Theorem}[section]
\newtheorem{prop}[theorem]{Proposition}

\theoremstyle{definition}
\newtheorem{defin}[theorem]{Definition}
\newtheorem{example}[theorem]{Example}

\newcommand{\ts}{\hspace{0.5pt}}
\newcommand{\nts}{\hspace{-0.5pt}}
\newcommand{\RR}{\mathbb{R}\ts}

\newcommand{\ZZ}{\mathbb{Z}}
\newcommand{\NN}{\mathbb{N}}
\newcommand{\QQ}{\mathbb{Q}}

\newcommand{\cA}{\mathcal{A}}

\newcommand{\cF}{\mathcal{F}}
\newcommand{\cL}{\mathcal{L}}
\newcommand{\cO}{\mathcal{O}}

\newcommand{\vL}{\varLambda}

\newcommand{\oplam}{\mbox{\Large $\curlywedge$}}

\newcommand{\exend}{\hfill$\Diamond$}

\newcommand{\myfrac}[2]{\frac{\raisebox{-2pt}{$#1$}}
	{\raisebox{0.5pt}{$#2$}}}

\numberwithin{equation}{section}

\makeatletter
\renewcommand{\@captionfont}{\small}
\makeatother

\DeclareMathOperator{\dotcup}{\dot\cup}

\newcommand{\defeq}{\mathrel{\mathop:}=}

\begin{document}
	
\title[A naturally appearing family of Cantorvals]{A naturally
  appearing family of Cantorvals}
	
\author{Michael Baake}
\address{Fakult\"at f\"ur Mathematik, Universit\"at Bielefeld,
	\newline \indent Postfach 100131, 33501 Bielefeld, Germany}
\email{$\{$mbaake, jmazac$\}$@math.uni-bielefeld.de}
	
\author{Anton Gorodetski}
\address{Department of Mathematics, University of California,
	\newline \indent Rowland Hall, Irvine, CA 92697, USA}
\email{asgor@uci.edu}
	
\author{Jan Maz\'{a}\v{c}}



\begin{abstract}
  The aim of this note is to show the existence of a large family of
  Cantorvals arising in the projection description of primitive
  two-letter substitutions. This provides a common and naturally
  occurring class of Cantorvals.
\end{abstract}

\maketitle
\section{Introduction}

A Cantorval is a subset of the real line enjoying properties both of
an interval and a~Cantor set. Formally, a Cantorval is a compact
subset of a line with dense interior such that none of its connected
components is isolated.  The name Cantorval was introduced by Mendes
and Oliviera \cite{MenOli} while studying sums of two homogeneous
Cantor sets. In such a scenario, several phenomena may occur.  One can
end up with a~Cantor set or with an interval --- or one obtains a
Cantorval.

Another example is given in \cite{MorMor}, where the authors show that
the sum of two dynamically defined Cantor sets with the sum of their
Hausdorff dimensions slightly larger than $1$ will typically be a
Cantorval. Such Cantorvals also appear in the study of achievement
sets. Here, an achievement set of a given vanishing sequence is the
set of the sums of its summable subsequences. The exploration of such
objects was initiated by Kakeya \cite{Kakeya}; see Nitecki
\cite{Nitecki} for a survey and G\l\c{a}b and Marchwicki \cite{GlMar}
for details. For a specific subclass of summable sequences, the
achievement set can be described as an attractor of an affine iterated
function system \cite{BFGPS}. Its topological and measure-theoretical
properties were studied in \cite{BBFS}.

It is conjectured that Cantorvals may appear as the spectrum of a
discrete two-dimensional Schr\"{o}dinger operator with separable
aperiodic potential given by the Fibonacci sequence in both directions
in an intermediate coupling regime. For weak coupling, it is known
that the spectrum is an interval, while it becomes a zero-measure
Cantor set in the large coupling regime. The intermediate remains
unclear; see \cite{DGS15,DEG15} for the known results and for images
of the spectra.

The aim of this short note is to present a new family of Cantorvals
which \emph{naturally} arise from studying geometric, self-similar
realisations of aperiodic sequences (one-dimensional quasicrystals)
with two symbols and their description as regular model sets. We
briefly introduce the projection method and present simple examples of
Cantorvals.  Finally, we present a theorem with conditions for the
occurrence of Cantorval in this context.

If the boundary of a Cantorval is a Cantor set, it is sometimes
referred to as an $M$-Cantorval (or, by some other authors, as
a \emph{symmetric} Cantorval). In what follows, we always mean
$M$-Cantorvals when we speak about Cantorvals.

\section{Primitive substitutions and the projection formalism}

Here, we are concerned with binary substitution rules, with the
alphabet $\cA = \{ a,b \}$.  For general background and details, we
refer to \cite[Sec.~4]{TAO} and references therein.

\begin{example}[Fibonacci substitution]
Consider
\[
  \varrho_{_\mathrm{F}} \, :\, \biggl\{ \!\begin{array}{l}
		a \mapsto ab, \\ b\mapsto a \ts , \end{array}
\]
which we abbreviate as
$\varrho_{_\mathrm{F}}= (ab,a) = \bigl( \varrho_{_\mathrm{F}} (a),
\varrho_{_\mathrm{F}} (b) \bigr)$ from now on.  One can iterate the
mapping
\[
  a\, \xrightarrow{\varrho_{_\mathrm{F}}}\,
  \varrho_{_\mathrm{F}}(a)=ab \, \xrightarrow{\varrho_{_\mathrm{F}}}\,
  \varrho_{_\mathrm{F}}(a)\varrho_{_\mathrm{F}}(b) = aba
  \,\xrightarrow{\varrho_{_\mathrm{F}}} \cdots
\]
and obtain a one-sided fixed point, that is, $\omega \in \cA^{\NN_0}$
such that $\varrho_{_\mathrm{F}}(\omega) = \omega$. Moreover, we can
extend it to a bi-infinite sequence by taking the legal seed $a|a$ as
a starting point (with $|$~denoting the position of zero) and
iterating
\[
  a|a \, \xrightarrow{\varrho_{_\mathrm{F}}} \, ab|ab\,
  \xrightarrow{\varrho_{_\mathrm{F}}}\, aba|aba\,
  \xrightarrow{\varrho_{_\mathrm{F}}}\, abaab|abaab\,
  \xrightarrow{\varrho_{_\mathrm{F}}}\cdots
  \xrightarrow{\varrho_{_\mathrm{F}}}\, w\,
  \xrightarrow{\varrho_{_\mathrm{F}}} \, \varrho_{_\mathrm{F}}(w)\,
  \xrightarrow{\varrho_{_\mathrm{F}}}\, w\,
  \xrightarrow{\varrho_{_\mathrm{F}}} \cdots
\]
which converges towards a $2$-cycle, each element of which is a fixed
point under the square of the substitution,
$\varrho^{2}_{_\mathrm{F}} = (aba, ab)$. The two elements of the
$2$-cycle differ only in the first two positions left of the marker
(where one sees either $ab$ or $ba$) and agree everywhere else. They
are an example of a proximal (in fact, asymptotic) pair.  \exend
\end{example}

Every binary substitution $\varrho$ induces a non-negative, integer
matrix $M_{\varrho} = (M_{ij})^{}_{1\leqslant i,j \leqslant 2}$ via
\[
  M_{ij} \, = \,
  \mbox{number of occurrences of letter } i \mbox{ in } \varrho(j).
\]
The substitution matrix gives insight into the nature of a given
substitution, for instance via some standard \emph{Perron--Frobenius}
(PF) arguments, some of which will also appear below.

\begin{defin}
  A substitution $\varrho$ is called \emph{primitive}, if
  $M_{\varrho}$ is primitive (so all entries of some power of
  $M_{\varrho}$ are positive). The substitution is \emph{unimodular},
  if $\det ( M_{\varrho} ) = \pm 1$.
\end{defin}

From now on, we always assume that our substitution is primitive.

\begin{example}\label{ex:more}
  The substitution matrix of $\varrho_{_\mathrm{F}}$ reads
  $\bigl(\begin{smallmatrix} 1 & 1 \\ 1& 0
  \end{smallmatrix}\bigr)$ and its PF eigenvalue is
  the golden ratio, $\tau = \tfrac{1+\sqrt{5}}{2}$. For the
  substitution
\[
    \varrho^{}_{2} \, = \, ( bba, ab) \ts ,
\]
the matrix reads $\bigl(\begin{smallmatrix} 1 & 1 \\ 2 & 1
\end{smallmatrix}\bigr)$, with PF eigenvalue $1+\sqrt{2}$. In
both cases, the substitution is unimodular, so
$\lvert\det(M_{\varrho})\rvert = 1$.  An example of a
\emph{non-unimodular} (nu) substitution is
\[
  \varrho^{}_{\mathrm{nu}} \, = \, ( aaba, aa)
\]
with
$M_{\varrho^{}_{\mathrm{nu}}} = \bigl(\begin{smallmatrix} 3 & 2 \\ 1 & 0
\end{smallmatrix}\bigr)$, $\det ( M_{\varrho^{}_{\mathrm{nu}}} ) = -2$
and $\lambda^{}_{\mathrm{PF}} = \frac{1}{2} ( 3 + \sqrt{17}\,)$.
\exend
\end{example}

A substitution matrix generally does not determine the
substitution. Indeed, consider
\begin{equation}\label{eq:Fibo-scrambled}
  \varrho^{2}_{_\mathrm{F}} \, =  \, ( aba, ab)  \quad \text{and} \quad
  \widetilde{\varrho} \, = \, ( aab,  ba )\ts ,
\end{equation}
which both share the substitution matrix
$\bigl(\begin{smallmatrix} 2 & 1 \\ 1& 1 \end{smallmatrix}\bigr)$.
Nevertheless, they are fundamentally different as we shall see shortly.

Every bi-infinite fixed point of a binary substitution can be
understood as a tiling of the real line --- it suffices to assign an
interval of a given length to each letter. The choice of the interval
lengths is not unambiguous, but there is (up to a scaling) only one
which turns the substitution $\varrho$ into a~geometric
\emph{inflation rule} (which will also be denoted by $\varrho$). The
latter first inflates all tiles by~$\lambda$, the PF eigenvalue of
$M_{\varrho}$, and then subdivides them according to the substitution
rule. It is not hard to see that the so-called \emph{natural tile
  lengths} are proportional to the entries of the left PF eigenvector
of $M_{\varrho}$. It is often a convenient choice for the shortest
tile to be of length 1.  The right PF eigenvector gives the relative
letter/tile frequencies (in statistical normalisation).

\begin{figure}[t]
  \includegraphics[width=0.6\textwidth,clip]{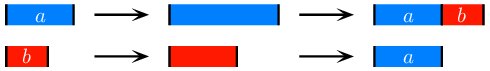}
  \caption{The geometric inflation rule for the Fibonacci tiling
    with natural tile sizes. \label{fig:Fibo} }
\end{figure}

In the case of $\varrho_{_\mathrm{F}}$, the letters $a$ and $b$ can be
replaced by intervals of length $\tau$ and $1$, respectively, which
are our prototiles. This turns the Fibonacci substitution rule into
the self-similar inflation rule of Figure~\ref{fig:Fibo}, which we
also call $\varrho_{_\mathrm{F}}$ for simplicity.  By decorating every
tile with a unique control point (for example, taking the left
endpoint of every interval), the resulting tiling can be described as
a point set $\vL$ (to be more precise, as a~disjoint union of two
point sets, so $\vL = \vL_a \, \dot{\cup}\, \vL_b$ in obvious
notation). Under some mild conditions, these points sets can be
described via the \emph{cut-and-project method} as follows.

First, consider the \emph{return vectors}, the translations that shift
a tile (or a patch of tiles) to another occurrence within the same
tiling, and consider all their integer combinations. The resulting set
is a finitely generated $\ZZ\ts$-module, known as the \emph{return
  module}.  In the case of a binary alphabet, its rank is at most
two. If our sequence is non-periodic, the rank equals two. Due to the
self-similar nature, the same module will be generated by the return
vectors for larger patches, because \emph{any} patch occurs within
$\varrho^n (a)$ for some $n\in\NN$ due to primitivity, which shows its
universal property.

Explicitly, if one follows the above convention and sets the tile
lengths to $\beta$ and $1$ (again induced by the entries of the left
PF eigenvector), the return module is
$\ZZ[\beta] \subset \QQ(\lambda)$, where $\lambda$ is the PF
eigenvalue of the given substitution. Note that $\beta$ may equal
$\lambda$, as in our guiding Fibonacci example, but this need not be
true in general, as for $\varrho^{}_{2}$ from
Example~\ref{ex:more}. The set $\ZZ[\beta]$ is a dense subset of
$\RR$, which can be described as the projection of a
higher-dimensional lattice, the \emph{Minkowski embedding} of
$\ZZ[\beta]$; see \cite[Sec.~3.4]{TAO}.

In our case, $\lambda$ is a real quadratic irrationality, so there are
two embeddings of $\QQ(\lambda)$ into the reals (given by the identity
and by the non-trivial Galois isomorphism that sends $\lambda$ to its
algebraic conjugate). Moreover, if the substitution is unimodular, the
scaling (or inflation) factor $\lambda$ is a unit, and the Minkowski
embedding of the return module $\ZZ[\beta]$ is indeed a lattice in
$\RR^2$; see \cite{TAO,Bernd} for background and details.

This gives a rise to a \emph{Euclidean cut-and-project scheme} (CPS)
as follows,
\begin{equation}\label{eq:CPS}
  \renewcommand{\arraystretch}{1.2}
  \begin{array}{ccccc@{}l}
   \RR & \xleftarrow{\;\;\; \pi \;\;\; }
       & \RR \nts\nts \times \nts\nts \RR &
       \xrightarrow{\;\: \pi^{}_{\text{int}} \;\: } & \RR & \\
       \cup & & \cup & & \cup  & \hspace*{-2ex}
       \raisebox{1pt}{\text{\scriptsize dense}} \\
       \pi (\cL) & \xleftarrow{\;\ts 1-1 \;\ts } & \cL &
       \xrightarrow{ \qquad } &\pi^{}_{\text{int}} (\cL) & \\
       \| & & & & \| & \\ L = \ZZ[\beta] &
       \multicolumn{3}{c}{\xrightarrow{\qquad\quad\quad \  \star
                     \quad\qquad\qquad}}
        &  {L_{}}^{\star\nts} = \ZZ[\beta^{\star}]  &  \end{array}.
   \renewcommand{\arraystretch}{1}
\end{equation}
Here, $\cL$ is the Minkowski embedding of $\ZZ[\beta]$ as mentioned
above, $\pi$ is the natural projection to the first coordinate (in
direct or physical space) such that $\pi\bigl|_{\cL}$ is injective,
and $\pi^{}_{\text{int}}$ is the projection to the second coordinate
(in internal space) where we assume that $\pi^{}_{\text{int}}(\cL)$ is
dense in $\RR$. Under these conditions, the star map
$\star : \QQ (\lambda) \longrightarrow \RR$ is well defined and given
by the non-trivial Galois automorphism of $\QQ(\lambda)$.

For the Fibonacci tiling, the return module is
$\ZZ[\tau] = \{m+n\tau : m,n\in\ZZ \} = \ZZ \oplus \tau\ZZ$, and the
Minkowski embedding of $\ZZ[\tau]$ into $\RR^2$ reads
\[
  \big\{ (x,x^{\star}) : x\in \ZZ[\tau] \big\} \, = \,
  \ZZ (1,1) \oplus \ZZ(\tau, 1-\tau)
\]
with $(a+b\tau)^{\star} = (a {+} b) - b\ts \tau$ for all $a,b\in \QQ$.

Having a CPS of the form \eqref{eq:CPS}, one can construct a
one-dimensional \emph{cut-and-project set} $\oplam (W)$ using a
suitable bounded subset $W$ --- the \emph{window} --- in internal
space via
\[
  \oplam (W) \, \defeq \,
  \{x\in L \, : \, x^{\star} \in W \}
  \, \subset \, \ZZ[\beta].
\]
So, $\oplam (W)$ consists of the $\pi$-projections of all lattice
points in the strip $\RR{\times} W$. If $W$ is a~relatively compact
set with non-empty interior, $\oplam (W)$, as well as $\oplam (W) +t$
with $t\in \RR$, is called a \emph{model set}. The latter is called
\emph{regular} when $\partial W$ has zero Lebesgue measure in internal
space. Note that the properties of $W$ are, not surprisingly,
reflected in the properties of $\oplam (W)$. More details and further
references can be found in \cite[Sec.~7]{TAO} or \cite{Moo97}.

\section{Pisot substitutions}

A natural question is how to determine a window $W$ such that
$\oplam (W)$ coincides with the set of control points $\vL$ (assuming
$0\in\vL$ without loss of generality) of a tiling induced by a given
substitution $\varrho$. The answer is not yet known in full
generality. However, in the 2-letter case under one additional
assumption on the scaling factor $\lambda$ (the PF eigenvalue of
$M_{\varrho}$), there exists a way to construct the window $W$ such
that $\vL = \oplam (W)$, or that $\vL \subseteq \oplam (W)$ with
$\oplam (W) \setminus \vL$ being a point set of density $0$.

The extra ingredient that allows us to proceed further is that
$\lambda$ is an algebraic integer with a particular property.  Recall
that $\lambda \in \RR$ is a \emph{Pisot--Vijayaraghavan number}, or PV
number for short, if it is an algebraic integer larger than $1$ such
that all its algebraic conjugates have modulus less than $1$, that is,
lie inside the open unit disk. If the PF eigenvalue attached to a
substitution $\varrho$ is a PV number, $\varrho$ is called a
\emph{Pisot substitution}. While it is still open (but widely believed
to be true in some form) whether every primitive Pisot substitution
with irreducible characteristic polynomial of its substitution matrix
gives rise to a regular model set, this is known to hold for binary
Pisot substitutions \cite{HS}.

For the binary case, a concrete embedding construction gives the
description of the self-similar inflation tiling as a model set for the
left endpoints of the tiles, as we first explain for the Fibonacci
tiling.

\begin{example}
  Consider the Fibonacci tiling, and let $\vL_{i}$ denote the set of
  all left endpoints of tiles of type $i \in\cA=\{a,b\}$, so
  $\vL = \vL_a \dotcup \vL_b$. Since the set $\vL$ is invariant under
  the action of the inflation, we get
\[
  \vL_a \,  = \, \tau \vL_a \,\dotcup \, \tau\vL_b
  \quad \text{and} \quad
  \vL_b \,  = \, \tau \vL_a + \tau  \ts ,
\]
where we momentarily suppress the possible objection that we should
rather work with $\varrho^2_{_\mathrm{F}}$.  The sets on the RHS
encode the positions of tiles in the inflated version of the
tiling. For example, the second equation states that all points of
type $b$ are obtained by taking all positions of $a$-points, inflating
them by $\tau$ and moving them by $\tau$ as the inflation rule in
Figure~\ref{fig:Fibo} suggests.  Taking the $\star$-image of these
equations and the closure of the sets, denoted by
$W_{i} = \overline{\vL^{\star}_{i}}$, yields
\begin{align*}
	W_a & = -\myfrac{1}{\tau} W_a \,\cup \, -\myfrac{1}{\tau}W_b, \\
	W_b & = -\myfrac{1}{\tau} W_a -\myfrac{1}{\tau},
\end{align*}
where we used that $\tau^{\star} = 1 - \tau = -\tfrac{1}{\tau}$.  This
constitutes an \emph{iterated function system} (IFS) for a pair of
non-empty, compact subsets of $\RR$, that is, an IFS on the space
$(\mathcal{K}\RR)^2$ equipped with the Hausdorff metric \cite{Wicks}.
We need the Pisot property of $\tau$, where
$|\tau^{\star}| = \tau-1 \approx 0.618\dots$, to see that the system
is contractive in the Hausdorff metric. Therefore, due to Hutchinson's
theorem \cite{Hutch81} and its extension to this setting
\cite{MauWil}, which is a version of Banach's contraction principle,
the system has the unique fixed point
\[
  W_a \, = \, [\tau-2,\tau -1], \qquad  W_b \, = \, [-1,\tau-2] \ts ,
\]
which can be checked by a direct calculation. The two compact sets are
closed intervals.

We can now resolve the subtlety mentioned above, which shows up in the
interpretation of the fixed point equation. Since our Fibonacci tiling
is one member of a singular $2$-cycle, the closed windows actually
describe the union of the two members. This adds one point in the
projection picture, which is not relevant to our discussion; see
\cite{BG-Rauzy} for details.
\exend
\end{example}

More generally, consider a binary Pisot substitution with inflation
factor $\lambda$, and turn to the corresponding self-similar inflation
rule with two prototiles (intervals) of natural length.  The above
procedure can be formulated via the set-valued \emph{displacement
  matrix} $T_{\varrho} = (T_{ij})^{}_{i,j\in \{a,b\}}$ with entries
\[
  T_{ij} \, = \, \{\mbox{all relative positions of a tile of type
       $i$ in } \varrho(j) \} \ts ,
\]
where $\varrho(j)$ refers to the level-$1$ supertile of type $j$ with
its decomposition into level-$0$ tiles. For instance, this gives
$T_{\varrho_{_\mathrm{F}}} = \left(\begin{smallmatrix} \{ 0 \} & \{ 0
    \} \\ \{ \tau \} & \varnothing \end{smallmatrix} \right) $ for the
Fibonacci inflation.  In general, we clearly have $|T_{ij}| =
M_{ij}$. Further, for $i\in \{a,b\}$, we get
\[
  \vL_i \, = \bigcup_{j\in\cA} \, \bigcup_{t \in T^{}_{ij}}
  \lambda \vL_j + t \ts ,
\]
and thus also
\begin{equation}\label{eq:IFS_gen}
  W_i \, = \bigcup_{j\in\cA} \, \bigcup_{t \in T^{}_{ij}}
  \lambda^{\star} W_j + t^{\star} \nts .
\end{equation}
This equation shows that the $W_i$ are dynamically defined (as
solutions of a graph-directed IFS, or GIFS). Usually, they are called
\emph{Rauzy fractals}; see \cite{PyFo} and references therein for
background.

Let us recall some known results about the $W_i$, which follow from
their definition as being the fixed point of an IFS; see
\cite[Cor.~6.66]{Bernd} for details and proofs, with plenty of further
references, and \cite{LGJJ93,FFIW} for important early examples.

\begin{prop}\label{prop:properties}
  Let\/ $W_a$ and\/ $W_b$ be two compact subsets of\/ $\RR$ that are
  the fixed point of the internal space IFS of a binary Pisot
  inflation rule. Then, the following properties hold.
\begin{enumerate}
\item $W_a$ and\/ $W_b$ have positive Lebesgue measure.
\item The unions on the RHS of Eq.~\eqref{eq:IFS_gen} are measure
  disjoint.
\item One has\/
  $\mathrm{int}(W_a) \ts \cap \ts \mathrm{int}(W_b) \, = \, \varnothing$.
\item The boundaries\/ $\partial W_a$ and\/ $\partial W_b$ have zero
  Lebesgue measure.
\item $W_a$ and\/ $W_b$ are perfect and topologically regular sets, that
  is, they are the closures of their interiors and do not possess
  any isolated points. \qed
\end{enumerate}
\end{prop}

The above construction in conjunction with
Proposition~\ref{prop:properties} shows that there are three different
kinds of resulting sets. Namely, each $W_i$ can be an interval, a
finite union of intervals, or the closure of an infinite (but
countable) union of disjoint intervals.  In the last case, the
boundary will have positive Hausdorff dimension, which forces the
windows to be Cantorvals.

Let us first consider the substitutions where the windows (for the
induced inflation rule) are simply intervals. This case has been
studied extensively, compare \cite{BEIR07,BFS12, Cant03,Lamb98}, where
the following criterion was derived.

\begin{prop}
  Under the assumptions of Proposition~$\ref{prop:properties}$, the
  two windows\/ $W_a$ and\/ $W_b$, as well as their union\/
  $W = W_a \cup W_b$, are connected sets, hence intervals, if and only
  if the underlying substitution\/ $\varrho$ is invertible, in the
  sense that\/ $\varrho$ is an automorphism of the free group\/
  $\cF^{}_2$ generated by\/ $\{a,b\}$.  \qed
\end{prop}

Let us illustrate this with our guiding example,
$\varrho_{_\mathrm{F}}$, where we know that the windows are
intervals. Thus, the Fibonacci substitution has to be
invertible. Indeed, one easily verifies that
$\varrho_{_\mathrm{F}}^{-1} = ( b, b^{-1}a)$ is the inverse of
$\varrho_{_\mathrm{F}}$, using that $b\ts b^{-1} = e$ is the neutral
element of $\cF^{}_2$.

Therefore, to get a Cantorval, one has to exclude all invertible
substitutions, because their fixed points are point sets which are
Delone (that is, they have bounded gaps and are relatively dense; see
\cite{TAO} for background) and give rise to windows that are
intervals. Likewise, one has to exclude all Delone sets that are
\emph{mutually locally derivable} (MLD; see \cite[Sec.~5.2]{TAO} for
details) to a fixed point set of an invertible inflation.  In short, a
point set $\vL_2$ is locally derivable from a point set $\vL_1$,
written as $\vL_1 \rightsquigarrow \vL_2$, when there is a rule of
fixed finite range that turns patches of $\vL_1$ into points of
$\vL_2$. The two sets are MLD when there are such finite-range rules
in both directions.

This concept is particularly powerful in the context of model sets.
Indeed, when comparing two model sets in the same CPS, say $\vL_1$ and
$\vL_2$ with corresponding windows $W_1$ and $W_2$, the derivability
$\vL_1 \rightsquigarrow \vL_2$ is equivalent to the property that
$W_2$ can be obtained from $W_1$ by \emph{finitely} many Boolean
operations with translates of $W_1$ (by projected lattice vectors),
which refers to intersections, unions, or complements; see
\cite{BSJ91, TAO} for more. So, any window occurring for a model set
that is MLD with a fixed point set of an invertible inflation rule can
at most be a \emph{finite} union of closed intervals.

\section{An example and a general result}

\begin{figure}[ht]
\centering
\includegraphics[width=0.91\textwidth]{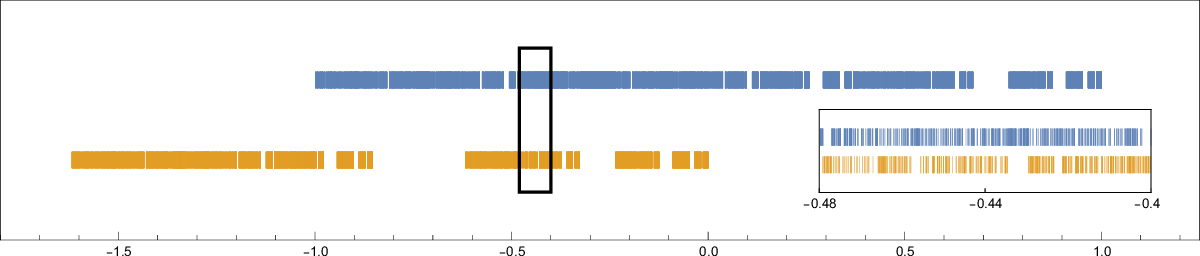}
\caption{\small The two windows (blue/top for control points of type
  $a$ and yellow/bottom for type $b$) for the tiling given by the
  substitution $\widetilde{\varrho}$ from
  Eq. \eqref{eq:Fibo-scrambled}. The inlay shows a stretched view of
  the marked region. We emphasise that the sets $W_a$ and $W_b$ are
  measure-theoretically disjoint, due to
  Proposition~\ref{prop:properties}, even though this is impossible to
  illustrate with the finite resolution of the plot. The image was
  produced via 10{\ts}000 steps of the random iteration algorithm (or
  `chaos game') for the attractor of an IFS, which relies on Elton's
  ergodic theorem \cite{Elton}. \label{fig:frac-window}}
\end{figure}

Let us discuss a binary Pisot inflation whose window is a Cantorval,
namely the substitution $\widetilde{\varrho} = (aab,ba)$ from
Eq.~\eqref{eq:Fibo-scrambled}. It shares the substitution matrix with
$\varrho_{_\mathrm{F}}^2$, so its PF eigenvalue is $\tau^2$, which is
a PV unit (since $(\tau^2)^{\star} = \tau^{-2} \approx 0.381966$).
The natural tile lengths are the same as for the Fibonacci example
(because $M^{}_{\varrho_{_\mathrm{F}}}$ and
$M^{}_{\widetilde{\varrho}} = M^{2}_{\varrho_{_\mathrm{F}}}$ have the
same eigenvectors) and so is the return module with its Minkowski
embedding. The displacement matrix becomes
$T_{\widetilde{\varrho}} = \left( \begin{smallmatrix} \{ 0, \tau \} &
    \{ 1 \} \\ \{ 2 \tau \} & \{ 0 \} \end{smallmatrix} \right)$, and
the iterated function system for the sets $W_a$ and $W_b$ now reads
\begin{equation}\label{eq:IFS-scrambled}
  \begin{split}
    W_a \, & = \, \tau^{-2} W_a \, \cup \,
        \tau^{-2} ( W_a - \tau ) \, \cup \,
        \tau^{-2} ( W_b + \tau^2 ) \ts , \\[1mm]
   W_b \, & = \, \tau^{-2} ( W_a - 2 \tau )
         \,\cup \, \tau^{-2}W_b \ts .
\end{split}
\end{equation}
The resulting windows are illustrated in Figure~\ref{fig:frac-window};
further details will appear in \cite{Jan}.

The picture reveals the Cantor-like structure which immediately
follows from the IFS~\eqref{eq:IFS-scrambled}, as the sets are
dynamically defined. The desired topological properties follow from
Proposition~\ref{prop:properties}. Moreover, one can explicitly show
the Cantor-like structure; compare \cite[Sec.~4.1]{LGJJ93}. A proof of
the fact that $W$ has infinitely many connected components is implicit
in \cite{Cant03}. One simple argument can be based on local
derivability as follows.

The Delone sets $\vL_a$ and $\vL$ are MLD, written as
$\vL_a \leftrightsquigarrow \vL$. Indeed, one has
$\vL \rightsquigarrow \vL_a$ from the observation that $\vL_a$ is the
subset of points of $\vL$ with a right neighbour at distance $\tau$,
and $\vL_a \rightsquigarrow \vL$ holds because the distance of an
$a$-point to its right neighbour is always $\tau$, $\tau +1$ or
$\tau+2$, which corresponds to $aa$, $aba$, $abba$, respectively. By
a~completely analogous argument, one has
$\vL_b \leftrightsquigarrow \vL$, and hence also
$\vL_a \leftrightsquigarrow \vL_b$ by transitivity.  Due to the MLD
criterion in terms of the windows mentioned earlier, see \cite{BSJ91}
as well as \cite[Rem.~7.6]{TAO}, the sets $W_a$, $W_b$ and $W$ must
then be of the same type, and their boundaries must have the same
Hausdorff dimensions. In particular, the dimensions are all zero or
all positive.

It is thus an important step to determine the Hausdorff dimension of
the boundary. If the windows were connected or a finite union of
connected components, the boundary dimension would (trivially) be
zero, and this is also true if the boundary would consist of countably
many points. In general, deriving the Hausdorff dimension of the
boundary may be technically difficult, but since the sets satisfy
Eq.~\eqref{eq:IFS-scrambled}, one can find a GIFS for
the boundary as well.

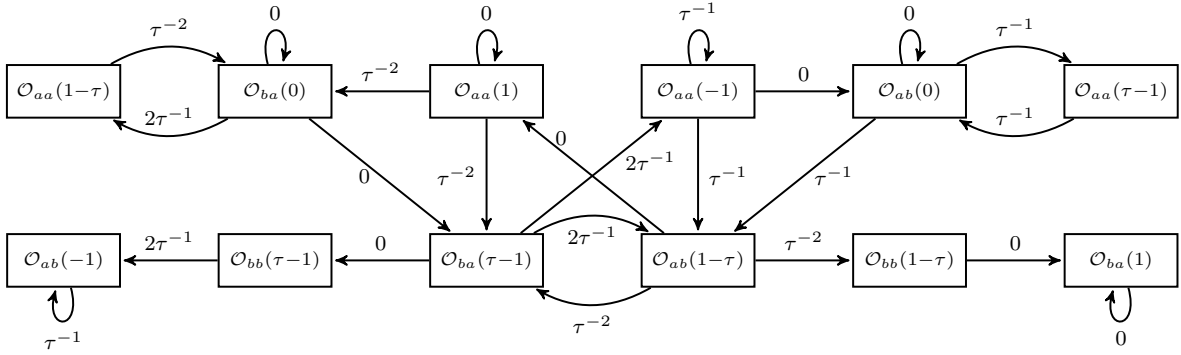
\begin{figure}[t]
\tikzset{every picture/.style={line width=0.75pt}} 
\begin{tikzpicture}[x=32pt,y=32pt,->,>=stealth',shorten >=1pt,
    auto,main node/.style={draw, rectangle,minimum width=1.5cm,
    minimum height=0.7cm, align=center,font=\tiny}]
 \node[main node] (a) at (0,0) {$\cO_{ab}(-1)$};
 \node[main node] (b) at (2.5,0) {$\cO_{bb}(\tau{-}1)$};
 \node[main node] (c) at (5,0) {$\cO_{ba}(\tau{-}1)$};
 \node[main node] (d) at (12.5,0) {$\cO_{ba}(1)$};
 \node[main node] (e) at (10,0) {$\cO_{bb}(1{-}\tau)$};
 \node[main node] (f) at (7.5,0) {$\cO_{ab}(1{-}\tau)$};
 \node[main node] (g) at (5,2) {$\cO_{aa}(1)$};
 \node[main node] (h) at (7.5,2) {$\cO_{aa}(-1)$};
 \node[main node] (i) at (2.5,2) {$\cO_{ba}(0)$};
 \node[main node] (j) at (10,2) {$\cO_{ab}(0)$};
 \node[main node] (k) at (0,2) {$\cO_{aa}(1{-}\tau)$};
 \node[main node] (l) at (12.5,2) {$\cO_{aa}(\tau{-}1)$};
		
\path
(a) edge [loop below] node {\scriptsize $\tau^{-1}$} (a)
(b) edge node [above] {\scriptsize $2\tau^{-1}$} (a)
(c) edge node [above] {\scriptsize $0$} (b)
(c) edge [bend left] node [below ] {\scriptsize $2\tau^{-1}$} (f)
(c) edge  node [right,xshift=0.3cm,yshift=0.2cm]
         {\scriptsize $2\tau^{-1}$} (h)
(d) edge [loop below] node {\scriptsize $0$} (d)
(e) edge node [above] {\scriptsize $0$} (d)
(f) edge node [above] {\scriptsize $\tau^{-2}$} (e)
(f) edge [bend left] node [below] {\scriptsize $\tau^{-2}$} (c)
(f) edge  node [left,yshift=0.5cm,xshift=-0.2cm] {\scriptsize $0$} (g)
(g) edge [loop above] node [above] {\scriptsize $0$} (g)
(g) edge  node [left] {\scriptsize $\tau^{-2}$} (c)
(g) edge  node [above] {\scriptsize $\tau^{-2}$} (i)
(h) edge [loop above] node [above] {\scriptsize $\tau^{-1}$} (h)
(h) edge  node [right] {\scriptsize $\tau^{-1}$} (f)
(h) edge  node [above] {\scriptsize $0$} (j)
(i) edge [loop above] node [above] {\scriptsize $0$} (i)
(i) edge [bend left]  node [above] {\scriptsize $2\tau^{-1}$} (k)
(i) edge  node [left] {\scriptsize $0$} (c)
(j) edge [loop above] node [above] {\scriptsize $0$} (j)
(j) edge [bend left]  node [above] {\scriptsize $\tau^{-1}$} (l)
(j) edge  node [right] {\scriptsize $\tau^{-1}$} (f)
(k) edge [bend left]  node [above] {\scriptsize $\tau^{-2}$} (i)
(l) edge [bend left]  node [above] {\scriptsize $\tau^{-1}$} (j);
\end{tikzpicture}	
\caption{The GIFS for the boundary substitution induced by
  $\widetilde{\varrho}$. The edge labeling corresponds to translates
  in the corresponding equations, as detailed in the text; see also
  \cite{Bernd}. The GIFS derived in \cite{FFIW} does not provide the
  full information about the geometry, though it also allows to
  compute the boundary dimension. \label{fig:bd}}
\end{figure}

Here, we present the boundary graph for the inflation
$\widetilde{\varrho}$, following the approach of \cite{Bernd}, where
one can find further details and examples. Denote by
$\cO_{\alpha \beta}(x)$ the ordered intersection
$W_{\alpha} \cap (W_{\beta} + x^{\star})$, which is called
$\varXi(\alpha,\beta,x)$ in \cite{Bernd}.  Then, the boundaries are
given as
\[
  \partial W_{\alpha} \, = \bigcup_{\beta\in\cA}
  \bigcup_x \cO_{\alpha \beta}(x) \ts ,
\]
where the sets $\cO_{\alpha \beta}(x)$ satisfy the GIFS encoded in
Figure~\ref{fig:bd}; see \cite{Jan} for more. Here, $x$ runs through
the non-zero vectors of shortest length in the return module of the
dual tiling (in internal space), modulo some reduction that is
described in detail in \cite[p.~178 and Rem.~ 6.102 and 6.103]{Bernd}.

This system of equations can be reduced further, without affecting the
result for the Hausdorff dimension, by identifying trival ones (such
as loops), solve them, and insert the result to simplify others. Also,
some symmetries can be employed, and one finally obtains
\begin{align*}
  \cO_{ba}(\tau{-}1)\,
  & = \, \bigl\{-\tau^{-1} \bigr\} \, \cup \, \tau^{-2}
    \cO_{aa}(-1)-2\tau^{-1} \, \cup \, \tau^{-2}
    \cO_{ba}(\tau{-}1)-\tau^{-1},\\[1mm]
  \cO_{aa}(-1)\,
  & = \,\tau^{-2} \cO_{aa}(-1)-\tau^{-1} \, \cup \, \tau^{-2}
    \cO_{ba}(\tau{-}1) \, \cup \, \tau^{-2} \cO_{ab}(0), \\[1mm]
  \cO_{ab}(0)\,
  & = \,\tau^{-2} \cO_{ab}(0) \, \cup \, \tau^{-2}
    \cO_{aa}(\tau{-}1) -\tau^{-1} \, \cup \, \tau^{-2}
    \cO_{ba}(\tau{-}1), \\[1mm]
  \cO_{aa}(\tau{-}1)\, & = \,\tau^{-2} \cO_{ab}(0)-\tau^{-1},
\end{align*}
together with $ \cO_{ab}(-1)=\{-1\}$, $\cO_{bb}(\tau{-}1) = \{-\tau\}$
and the general symmetry relation
$\cO_{\alpha \beta}(x) = \cO_{ \beta \alpha}(-x) + x^{\star}$.  The
spectral radius of the adjacency matrix of the boundary graph in
Figure \ref{fig:bd} provides the scaling factor needed for the
dimension calculation.  In the case of~$\widetilde{\varrho}$, this
spectral radius is $1+\sqrt{2}$. Thus, the Hausdorff dimension of the
boundary (which is the same for $\partial W_a$, $\partial W_b$ and
$\partial W$ as explained above) is
\[
  d_{_\mathrm{H}} \, = \, \frac{\log (1+\sqrt{2}\, )}{2\log (\tau)}
  \, \approx \, 0,915{\ts} 785{\ts} 46\dots \ts ,
\]
which can already be found in \cite{LGJJ93,FFIW}. One can
alternatively employ the \emph{orbit separation dimension} (OSD) from
\cite{BGG}. This confirms the non-trivial structure of the window
itself, and one concludes that $W$ is a~Cantorval.

These ideas can be applied to any such substitution, including the
substitution $\varrho^{}_{2}$ mentioned earlier, and then give the
following general result.

\begin{theorem}
  Let\/ $\varrho$ be a primitive, unimodular Pisot inflation on a
  binary alphabet. Consider its realisation as a model set with
  window\/ $W\!$. If\/ $\partial W$ has positive Hausdorff dimension,
  $W$ is a~Cantorval.  \qed
\end{theorem}

Note that the inflation structure gives not only an IFS for the
windows, but also induces a GIFS for the \emph{boundaries} of them. If
the Hausdorff dimension is positive, this GIFS establishes the Cantor
structure of the boundaries, and hence the Cantorval property.

Many concrete examples can be constructed, always via the same
recipe. One starts with a binary Pisot substitution, then takes some
power of it (if needed), and flips letters. This will typically change
the LI (local indistinguishability) class, as well as the frequency
module, and often leads to windows with fractal boundary, then giving
Cantorvals.

When one goes beyond binary alphabets, but sticks to unimodular
substitutions, the windows will generally be Rauzy fractals in $\RR^m$
for $m>1$, and more phenomena are possible. The closest analogue of
the Cantorval structure is the case that a window can consist of
countably, but infinitely many disjoint pieces; see \cite{ST} for
some partial classification.

The situation gets more complicated when one goes to
higher-dimensional inflation tilings. Indeed, already when looking in
direct product variations or more general higher-dimensional inflation
tilings, a plethora of new types of Rauzy fractals will emerge. As far
as we are aware, no classification is in sight, though the OSD from
\cite{BGG} can certainly help to analyse such cases.

\section*{Acknowledgements}

It is our pleasure to thank Franz G\"{a}hler for valuable discussions,
and two anonymous referees for their helpful suggestions to improve
our presentation.  This work was supported by the German Research
Council (Deutsche Forschungsgemeinschaft, DFG) under SFB-TRR 358/1
(2023) -- 491392403 and by NSF grant DMS-2247966.  

\end{document}